\documentclass{article}

\usepackage{amsmath}
\usepackage{amsfonts}
\usepackage{amssymb}
\usepackage{mathrsfs}
\usepackage{color,graphicx}

\newcommand   \Integers {\mathbb Z}

\newcommand   \fieldChar {\mathtt{\,p\,}}
\newcommand   \finiteIntegerField {\basel{\mathbb Z}{\fieldChar}\,} 
\newcommand   \scalars  {\mathbb F}

\newcommand	  \rationals {\mathbb Q}

\renewcommand   \mod  {~\textsf{mod}~ }
\renewcommand   \gcd  {~\textsf{gcd}~ }
\newcommand   \lcm  {~\textsf{lcm}~ }

\newcommand  \singlevariablepolynomials[2] {#1{\mathbf{[}#2\mathbf{]}}}
\newcommand  \basel[2]{#1_{_{#2}}}

\newcommand  \tab  {\hspace*{0.5cm}}

\newcommand  \shiftright  {\hspace*{2.0cm}}

\newtheorem{theorem}{Theorem}

\newcommand \proof {\rm{\textbf{Proof.}}~~}

\newcommand {\qed} {\hfill{} $\square$}

\newcommand  \bglb {\big (}
\newcommand  \bgrb {\big )}

\begin{document}

\title{\textcolor{magenta}{\bf{Factorization of Polynomials over the Field of Rational Numbers}}}

\author{ \\
\textcolor{green}{\bf{Duggirala Meher Krishna}}\\
\textcolor{cyan}{\small{Gayatri Vidya Parishad College of Engineering (Autonomous)}} \\
\textcolor{cyan}{\small{Madhurawada, VISAKHAPATNAM -- 530 048, Andhra Pradesh, India}} \\
\textcolor{cyan}{\small{E-mail ~: \tab duggiralameherkrishna@gmail.com}}\\
 \\
 and \\
 \\
\textcolor{green}{\bf{Duggirala Ravi}}\\
\textcolor{cyan}{\small{Gayatri Vidya Parishad College of Engineering (Autonomous)}} \\
\textcolor{cyan}{\small{Madhurawada, VISAKHAPATNAM -- 530 048, Andhra Pradesh, India}} \\
 \textcolor{cyan}{\small{E-mail ~: \tab ravi@gvpce.ac.in; \tab \tab \tab duggirala.ravi@yahoo.com};} \\
\textcolor{cyan}{\small{\shiftright  duggirala.ravi@rediffmail.com; \tab drdravi2000@yahoo.com}} 
}

\date{}

\maketitle

\begin{abstract}
% Text of abstract
In this paper, a randomized algorithm for deciding the irreducibility of an irreducible polynomial and factoring a reducible polynomial over the field of rational numbers is presented. The main idea underlying the algorithm is based on conversion of a given polynomial into a polynomial with integer coefficients and reduction to $\mod \fieldChar$, for several large prime numbers $\fieldChar$, without applying a lifting method. 
\\

% keywords here, in the form: keyword \sep keyword
\\

% PACS codes here, in the form: \PACS code \sep code
% MSC2010 codes here, in the form:  MSC2010 code \sep code
%\noindent {\em{MSC2010}:}~~~~ 11T06,~~12F05,~~12F10,~~12E20,~~12Y05. 

\end{abstract}

\textcolor{blue}{
\section{Introduction}}
In this paper, a randomized algorithm for deciding the irreducibility of an irreducible polynomial and factoring a reducible polynomial over the field of rational numbers is presented. The existing literature is based on lattice reduction over integers or application of suitable lifting methods, after obtaining the irreducible factors $\mod \fieldChar$, for a randomly chosen prime number $\fieldChar$. Some of the methods discussed in the literature are found in the references section. For a large prime number $\fieldChar$, the probability that an irreducible polynomial over the field of rational numbers, after its conversion into a polynomial with integer coefficients, remains irreducible when $\mod \fieldChar$ restriction is applied happens to be high, under the assumption of uniform probability over all possible monic irreducible polynomials of the degree at most that of the given polynomial. Thus, appropriate scalar multiples of the monic irreducible polynomials $\mod \fieldChar$ can be tested for being factors of the given polynomial with rational number coefficients over the field of rational numbers. The appropriate scalar to be used for  testing the divisibility is the coefficient of the leading degree term of the given polynomial, after its conversion to polynomial with integer coefficients. 

%\vspace*{-0.2cm}
 
\textcolor{blue}{ 
\section{\label{Sec-Formulation}The Polynomial Factorization Problem}}
Let $\Integers$ be the ring of integers,  $\fieldChar$ be a suitable large prime number, and $\finiteIntegerField$ be the finite field of integers $\mod \fieldChar$. As a special case, monic polynomials with integer coefficients are discussed, and the modifications required to apply for polynomials with rational coefficients are described later.

%\vspace*{-0.2cm}
\textcolor{blue}{\textbf{
\subsection{\label{SpecialCase}Monic Polynomials with Integer Coefficients}}}

The following result provides an estimate for the probability that a monic irreducible polynomial with integer coefficients remains irreducible, when $\mod \fieldChar$ restriction is applied.   
 \begin{theorem}
 \label{Theorem-1} Let $g(x)$ be a nonconstant monic irreducible polynomial of degree $s \geq 1$, with integer coefficients,  $\fieldChar$ be a suitable large prime number, and  $h(x) = g(x) \mod \fieldChar$. If $\gcd \bglb h(x),\, x \bgrb = 1$, $\deg\bglb h(x) \bgrb = s$   and every monic irreducible polynomial of degree at most $s$, other than a scalar multiple of $x$, is equally likely to be a factor of $h(x)$ over $\finiteIntegerField$,  then the probability that $h(x)$ remains an irreducible polynomial over $\finiteIntegerField$ is at least $ \bglb 1 - \frac{s+1}{\fieldChar-1} + {\mathcal O} \bglb \frac{s^{2}}{(\fieldChar-1)^{2}} \bgrb \bgrb$. 
 \end{theorem}
 \proof If $s = 1$, then there is nothing to prove. Let $s \geq 2$ and $\gcd \bglb h(x),\, x \bgrb = 1$. The irreducible factors of $h(x)$ over $\finiteIntegerField$ are those that are factors of the polynomials $(x^{\fieldChar^{i}-1}-1)$, for $1 \leq i \leq s$. Now, $h(x)$ is irreducible if and only if $\gcd\bglb h(x),  (x^{\fieldChar^{i}-1}-1)\bgrb = 1$, for $1 \leq i \leq s-1$, and $h(x) \, \mid \, (x^{\fieldChar^{s}-1}-1)$. The polynomial $L(x) = \prod_{i = 1}^{s-1}  (x^{\fieldChar^{i}-1}-1)$ is of degree $\sum_{i = 1}^{s-1} (\fieldChar^{i}-1) = \sum_{i = 0}^{s-1} \fieldChar^{i} - s  = \frac{\fieldChar^{s}-1}{\fieldChar-1} -s$.  The number of monic irreducible polynomials of degree at most $s-1$, other than scalar multiples of $x$, is at most $\frac{\fieldChar^{s}-1}{\fieldChar-1} -s$, since every such irreducible factor is a divisor of $L(x)$.  Let $R(x) = \gcd \bglb L(x), \, x^{\fieldChar^{s}-1} -1 \bgrb$ and $P(x) = \frac{x^{\fieldChar^{s}-1} -1}{R(x)}$. Then, the degree of $P(x)$ is at least $\fieldChar^{s}-1 - \bglb  \frac{\fieldChar^{s}-1}{\fieldChar-1} -s \bgrb$, and the number of monic irreducible polynomials of degree $s$ is at least  
 $\frac{\fieldChar^{s}-1 - \bglb  \frac{\fieldChar^{s}-1}{\fieldChar-1} -s \bgrb}{s}$. Now, the total number of monic irreducible polynomials, other than scalar multiples of $x$, is at most $\sum_{i = 1}^{s} \frac{1}{i} \times (\fieldChar^{i}-1) \leq \frac{(\fieldChar^{s}-1)}{s} +\sum_{i = 1}^{s-1} (\fieldChar^{i}-1)  = \frac{(\fieldChar^{s}-1)}{s} + \frac{\fieldChar^{s}-1}{\fieldChar-1} -s $. Thus, the probability that $h(x)$ is irreducible over $\finiteIntegerField$ is at least  $\frac{\fieldChar^{s}-1 - \bglb  \frac{\fieldChar^{s}-1}{\fieldChar-1} -s \bgrb}{s \cdot \bglb \frac{(\fieldChar^{s}-1)}{s} + \frac{\fieldChar^{s}-1}{\fieldChar-1} -s \bgrb}$
$ = $
$\frac{1 - \bglb  \frac{1}{\fieldChar-1} -\frac{s}{\fieldChar^{s}-1} \bgrb}{ 1 + s\cdot \bglb  \frac{1}{\fieldChar-1} -\frac{s}{\fieldChar^{s}-1} \bgrb}$  
 , which can be expanded into 
 $ \bglb 1 - \bglb  \frac{1}{\fieldChar-1} -\frac{s}{\fieldChar^{s}-1} \bgrb \bgrb \times \bglb 1 - s\cdot \bglb  \frac{1}{\fieldChar-1} -\frac{s}{\fieldChar^{s}-1} \bgrb + {\mathcal O}\bglb \frac{s^{2}}{(\fieldChar-1)^{2}}\bgrb \bgrb $, 
whereupon the conclusion follows.  \qed
\\

A more realistic estimate of the probability in the statement of the previous theorem is 
$\bglb \frac{\fieldChar^{s}-1 - \bglb  \frac{\fieldChar^{s}-1}{\fieldChar-1} -s \bgrb}{s\fieldChar^{s}} \bgrb = \frac{1}{s} \times \bglb 1 - \frac{1}{\fieldChar-1} + {\mathcal O} \bglb \frac{1}{(\fieldChar-1)^{2}} \bgrb \bgrb$, an estimate based on the assumption of uniform likelihood of occurrence of $h(x)$ as anyone of all possible monic polynomials of degree $s$. The converse of Theorem \ref{Theorem-1} also holds true: if $h(x)$ is an irreducible polynomial over $\finiteIntegerField$ of degree $s$, then so is also $g(x)$. The converse is useful for inferring the irreducibility of a given monic polynomial with integer coefficients, as will be discussed in next subsection. 
Let $f(x)$ be a nonconstant monic polynomial with integer coefficients of degree $d > 1$, with coefficient of the leading degree term equal to $1$. Let $f(x) = \prod_{i = 1}^{r} \basel{f}{i}(x)$, where $\basel{f}{i}(x)$ are monic irreducible factor polynomials with integer coefficients of degree $\basel{d}{i} \geq 1$, for $1 \leq i \leq r$, for some positive integer $r$.  
 The factoring algorithm tests for irreducible factors $\mod \fieldChar$, for several large prime integers $\fieldChar$, being factors of the given polynomial $f(x)$. By the previous theorem, the restrictions of the factor polynomials $\basel{f}{i}(x)$ to $\mod \fieldChar$, $1 \leq i \leq r$, remain irreducible over $\finiteIntegerField$, for suitably large prime numbers $\fieldChar$, with reasonably high probability.
 
%\vspace*{-0.2cm}
\textcolor{blue}{\textbf{
\subsection{\label{GeneralCase}Monic Polynomials with Rational Number Coefficients}}}
Let $f(x)$ be a nonconstant monic polynomial with rational coefficients of degree $d > 1$, with coefficient of the leading degree term equal to $1$. Let $f(x) = \prod_{i = 1}^{r} \basel{f^{\basel{\nu}{i}}}{i}(x)$, where $\basel{f}{i}(x)$ are distinct monic irreducible factor polynomials with rational number coefficients of degree $\basel{d}{i} \geq 1$, and $\basel{\nu}{i} \geq 1$ are the multiplicities, for $1 \leq i \leq r$, for some positive integer $r$. For some index $i$, where $1 \leq i \leq r$,  if $\basel{\nu}{i} > 1$, then $\basel{f^{\basel{\nu}{i}-1}}{i}(x)$ is a factor of $\gcd \bglb f(x), f'(x) \bgrb$, where $f'(x)$ is the formal derivative of $f(x)$, and so $\frac{f(x)}{\gcd \bglb f(x), f'(x) \bgrb} = \prod_{i = 1}^{r} \basel{f}{i}(x)$. The computation of $\gcd$ with the corresponding formal derivative can be repeated with $\gcd \bglb f(x), f'(x)\bgrb$ in place of $f(x)$, hopefully recovering some of the irreducible factors $\basel{f}{i}(x)$, for $1 \leq i \leq r$.

Let now  $f(x) = \prod_{i = 1}^{r} \basel{f}{i}(x)$, where $\basel{f}{i}(x)$ are distinct monic irreducible factor polynomials with rational number coefficients of degree $\basel{d}{i} \geq 1$, for $1 \leq i \leq r$, for some positive integer $r$. The rational number coefficients are assumed to be in the minimal representation of numerators and denominators, with no nontrivial integer factors between each pair.
 Then the $\lcm$ of denominators of the coefficients of $\basel{f}{i}(x)$ is also the coefficient  of the leading degree term of $\basel{f}{i}(x)$, after its conversion into a polynomial with integer coefficients with content $1$. Let $c$ be the $\lcm$ of the denominators of $f(x)$. The content being multiplicative,  $c\basel{f}{i}(x)$ is an irreducible polynomial of integer coefficients of degree $\basel{d}{i}$, $1 \leq i \leq r$. Then the factoring algorithm tests for the divisibility of $f(x)$ over the rational numbers, by the polynomial $h(x) =\bglb (c * g(x)) \mod \fieldChar \bgrb$, treating $h(x)$ as a polynomial with integer coefficients, for every monic irreducible factor $g(x)$ of $\bglb (c * f(x)) \mod \fieldChar \bgrb$ over $\finiteIntegerField$.

%\vspace*{-0.2cm}
\textcolor{blue}{
\section{\label{First-extension}Factorization over Algebraic Extensions of Rational Numbers}}
Let $\rationals$ be the field of rational numbers, $\singlevariablepolynomials{\rationals}{x}$ be the polynomial ring in the indeterminate $x$ with rational number coefficients, $\phi(˘x) \in \singlevariablepolynomials{\rationals}{x}$ be a monic irreducible polynomial of degree $k \geq 2$, and $\scalars = \singlevariablepolynomials{\rationals}{\alpha}/\phi(\alpha)$ be the extension field of dimension $k$ over $\rationals$ with arithmetic operations $\mod  \phi(\alpha)$.  The nonzero elements in $\scalars$ are assumed to be represented in the simplest form of their representation, as polynomials in $\alpha$ of degree at most $k-1$ with rational number coefficients.  The following algorithm can be used for testing the irreducibility and / or factorization, whichever is the case, of the polynomial $f(x)$ over $\scalars$:
\\

\begin{enumerate}
\item Let $\psi(x)$ be a polynomial with integer coefficients obtained by multiplying $\phi(x)$ with an appropriate integer. For several large prime numbers $\fieldChar$, such that $\psi(x)$ is irreducible over $\finiteIntegerField$, irreducibility of $f(x)$ over $\singlevariablepolynomials{\finiteIntegerField}{\gamma}/\psi(\gamma)$ is tested, and if $f(x)$ is found to be irreducible for any single large prime number  $\fieldChar$ in this step, the procedure is terminated.

\item If $f(x)$ is a polynomial with coefficients in $\rationals$, let $F(x) = f(x)$,   and the procedure moves to the next step, skipping the remaining part in this step. Now the symmetric product $F(x) = \prod_{i=0}^{k-1}\basel{\theta}{i}f(x)$ is computed,  where $\basel{\theta}{i}$, $0 \leq i \leq k-1$, are the distinct embeddings of $\scalars$ in the algebraic closure of $\rationals$. For elements in $\scalars$, the symmetric product as just discussed is the norm of an element, as found in field theory. The norm of an element in $\scalars$ is the determinant of the matrix corresponding to multiplication by the given element, with respect to any particular basis for $\scalars$ over $\rationals$. For extending the norm to polynomials over $\scalars$, a specific scheme is considered. Let ${\mathcal K}$ be the field of rational functions over $\rationals$, as formal expressions in a single indeterminate $x$, with $\rationals$ embedded in ${\mathcal K}$ by the identity morphism, and ${\mathcal L} = \singlevariablepolynomials{{\mathcal K}}{\alpha}/\phi(\alpha)$, assuming that $\phi(y)$ remains irreducible in $\singlevariablepolynomials{{\mathcal K}}{y}$. Then the simple extension $\scalars$ over $\rationals$ naturally extends to ${\mathcal L}$ over ${\mathcal K}$. The given polynomial $f(x)$ is suitably represented as an element in ${\mathcal L}$. Then the norm of $f(x)$ is computed for the extension ${\mathcal L}$ over ${\mathcal K}$. This method is very efficient, for computing the symmetric product. 

\item The polynomial $F(x)$ is a polynomial with rational number coefficients, and hence, can be factorized over $\rationals$ using the methods of the previous section.  Let $F(x) = \prod_{i = 1}^{r} \basel{G^{\basel{\nu}{i}}}{i}(x)$, where $\basel{G}{i}(x)$ are distinct irreducible factors of multiplicity $\basel{\nu}{i} \geq 1$, be the factorization of $F(x)$ into irreducible polynomials over $\rationals$. Now, let $\basel{f}{i}(x) = \gcd\bglb f(x), \basel{G}{i}(x) \bgrb$, so that $\basel{f}{i}(x)$ is the product of $\basel{\nu}{i}$ many irreducible factors that are related by the conjugacy of $\scalars$ over $\rationals$, $1 \leq i \leq r$.  For finding factors of multiplicity greater than $1$, the $\gcd$ of the given polynomial with its derivative is tested. For factoring a product $\basel{f}{i}(x)$ of $\basel{\nu}{i}$  many irreducible polynomials,  where $\basel{\nu}{i} > 1$, each of multiplicity $1$, but related among themselves by the conjugacy, the polynomial $\basel{f}{i}(x-\lambda \alpha)$ is tested for reducibility, for various values of $\lambda \in \rationals$, where $\alpha$ is a fictitious root of $\phi$. Some of the irreducible factors of $\basel{f}{i}(x-\lambda\alpha)$ may not be related to other irreducible factors by the conjugacy, and hence, can be factored by repeating the previous and current steps, for some particular $\lambda \in \rationals$. If $g(˘x)$ is a factor of $\basel{f}{i}(x-\lambda\alpha)$, then $g(x+\lambda \alpha)$ is a factor of $\basel{f}{i}(x)$, for $1 \leq i \leq r$.
\end{enumerate}
 
The algorithm can be used for finding the automorphisms by taking the symmetric product of $f(x)$, for various polynomials $f(x)$ with coefficients in $\scalars$, and factoring it over $\scalars$. This approach allows the solver to observe how the roots are related, allowing the automorphism group of $\scalars$ over $\rationals$ to be efficiently computed. It is also preferable to factorize $\phi(x)$ with respect to a fixed basis, such as $\{\alpha^{i}\,:\, 0 \leq i \leq k-1\}$, for finding the roots themselves in terms of any one fixed fictitious root $\alpha$ of $\phi(x)$.  Factorization over multiple algebraic extensions is recursively performed.

\textcolor{blue}{
 \section{Conclusions}}
 The aim of the paper is to provide a justification for a randomized polynomial factoring algorithm by considering restriction to $\mod \fieldChar$, for several suitably large prime numbers $\fieldChar$. A realistic estimate for the probability that an irreducible factor over integers of degree $s$ remains irreducible with $\mod \fieldChar$ is close to $\frac{1}{s}$, which is reasonably encouraging, for achieving success by repeated trials. Extension of the proposed method to factorization over algebraic extensions of the field of rational numbers is also possible. 
\\
\\

\end{document}